\documentclass[fleqn]{mat01}
\usepackage{times,mathtimy,amssymb,latexsym}
\begin{document}

\def\d{\mbox{\rm d}}
\def\e{\mbox{\rm e}}

\setcounter{page}{393} \firstpage{393}

\newtheorem{theore}{Theorem}
\renewcommand\thetheore{\arabic{section}.\arabic{theore}}
\newtheorem{theor}{\bf Theorem}
\newtheorem{prob}[theor]{\it Problem}
\newtheorem{exam}[theor]{Example}

\def\lemm{\trivlist \item[\hskip \labelsep{\it Lemma.}]}
\def\prop{\trivlist \item[\hskip \labelsep{\rm PROPOSITION.}]}

\def \diag  {\text {\rm diag}}
\def \cos   {\text {\rm cos}}
\def \Eig   {\text {\rm Eig}}
\def \Max   {\text {\rm Max}}
\def \Min   {\text {\rm Min}}
\def \min   {\text {\rm min}}
\def \Re    {\text {\rm Re}}
\def \dt    {\text {\rm dt}}
\def \diag  {\text {\rm diag}}
\def \lim   {\text {\rm lim}}
\def \sup   {\text {\rm sup}}
\def \tr    {\text {\rm tr}}
\def \wt    {\widetilde}
\def \qed   {\hfill \vrule height6pt width 6pt depth 0pt}

\title{Corners of normal matrices}

\markboth{Rajendra Bhatia and Man-Duen Choi}{Corners of normal
matrices}

\author{RAJENDRA BHATIA and MAN-DUEN CHOI$^{*}$}

\address{Theoretical Statistics and Mathematics Unit, Indian Statistical Institute, New~Delhi~110~016, India\\
\noindent $^{*}$Department of Mathematics, University of Toronto, Toronto M5S 2E4, Canada\\
\noindent E-mail: rbh@isid.ac.in; choi@math.toronto.edu\\[1.2pc]
\noindent {\it To Kalyan Sinha on his sixtieth
\vspace{-1pc}birthday}}

\volume{116}

\mon{November}

\parts{4}

\pubyear{2006}

\Date{}

\begin{abstract}
We study various conditions on matrices $B$ and $C$ under which they can be the
off-diagonal blocks of a partitioned normal matrix.
\end{abstract}


\maketitle


\noindent The structure of general normal matrices is far more
complicated than that of two special kinds~---~hermitian and
unitary. There are many interesting theorems for hermitian and
unitary matrices whose extensions to arbitrary normal matrices
have proved to be extremely recalcitrant (see e.g., \cite{1}). The
problem whose study we initiate in this note is another one of
this sort.

We consider normal matrices $N$ of size $2n,$ partitioned into blocks of size $n$ as
\begin{equation}
N = \begin{bmatrix} A & B \\[.2pc] C & D \end{bmatrix}.\label{eq1}
\end{equation}
Normality imposes some restrictions on the blocks. One such restriction is the
equality
\begin{equation}
\|B\|_{2} = \|C\|_{2} \label{eq2}
\end{equation}
between the {\em Hilbert--Schmidt (Frobenius) norms} of the off-diagonal blocks $B$
and $C$. If $T$ is any $m \times m$ matrix with entries $t_{ij},$ then
\begin{equation*}
\| T \|_2 = \left( \sum_{j=1}^m |t_{ij}|^2 \right)^{1/2}.
\end{equation*}
The equality \eqref{eq2} is a consequence of the fact that the Euclidean norm of the
$j$th column of a normal matrix is equal to the Euclidean norm of its $j$th row.

Replacing the Hilbert--Schmidt norm by another unitarily invariant norm, we may ask
whether the equality \eqref{eq2} is replaced by interesting inequalities. Let $s_1 (T)
\ge \cdots \ge s_m (T)$ be the singular values of $T$. Every unitarily invariant norm
$||| T |||$ is a symmetric gauge function of $\{ s_j (T) \}$ (see chapter IV of
\cite{1} for properties of such norms). Much of our concern in this note is with the
special norms
\begin{equation*}
\| T \|_2 = (\mbox{tr}\, T^{*} T)^{1/2} = \left (\sum_{j=1}^m s_j^2 (T) \right)^{1/2}
\end{equation*}
and
\begin{equation}
\| T \| = s_1 (T) = \sup_{x \in {{\Bbb C}^{m}}, \|x\|=1} \| Tx \|. \label{eq3}
\end{equation}
The latter is the norm of $T$ as a linear operator on the Euclidean space ${\Bbb C}^m.$ Clearly
\begin{equation}
\| T \| \leq \| T \|_2 \leq \sqrt{m} \| T \|, \label{eq4}
\end{equation}
for every $m \times m$ matrix $T.$

If the matrix $N$ in \eqref{eq1} is hermitian, then $C=B^*,$ and hence,
$|||C|||=|||B|||$ for all unitarily invariant norms. If $N$ is unitary, then
$AA^*+BB^*=A^*A+C^*C=I.$ Hence, the eigenvalues $\lambda_j$ satisfy the relations
\begin{align*}
\lambda_j (BB^*) &= \lambda_j (I-AA^*) = 1 - \lambda_j (AA^*)\\[.3pc]
&= 1 - \lambda_j (A^*A) = \lambda_j (I-A^*A) = \lambda_j (C^*C).
\end{align*}
Thus $B$ and $C$ have the same singular values, and again $|||B||| = |||C|||$ for all
unitarily invariant norms.

This equality of norms does not persist when we go to arbitrary normal matrices, as we
will soon see. From \eqref{eq2} and \eqref{eq4} we get a simple inequality
\begin{equation}
||B|| \leq \sqrt{n} \ ||C||.\label{eq5}
\end{equation}

One may ask whether the two sides of \eqref{eq5} can be equal, and that is the first
issue addressed in this note.

When $n=2,$ it is not too difficult to construct a normal matrix $N$ of the form
\eqref{eq1} in which $\|B\| = \sqrt{2} \|C\|.$ One example of such a matrix is
\begin{equation}
N = \begin{bmatrix} \frac{\begin{array}{rl}
 0 &0\\[.2pc]
 1 &0\end{array}}{\begin{array}{rl}
 0 &1\\[.2pc]
1 &0\end{array}} \vline \frac{\begin{array}{rl}
 \sqrt{2} &0\\[.2pc]
 0& 0\end{array}}{\begin{array}{@{ \ \ \ \ }cl}
0 &1\\[.2pc]
0 &0\end{array}}\end{bmatrix}.\label{eq6}
\end{equation}
\noindent When $n=3,$ examples seem harder to come by. One that
preserves some of the features of \eqref{eq6} is given by the
matrix
\begin{equation}
N = \begin{bmatrix} \frac{\begin{array}{lcr}  0 &
\sqrt{\frac{2}{\sqrt{3}}-1} & 0\\
& &\\
0 &0 &\sqrt{\frac{2}{\sqrt{3}}}\\
& &\\
\sqrt{\frac{2}{\sqrt{3}}+1} & 0 & 0\end{array}}{\begin{array}{lcccr}
 0 \quad \quad \: \: & &\quad \quad 0 && \quad \quad \: \:1\\
 &&&&\\
 &&&&\\
 0 & &\quad \quad 1& &0\\
 &&&&\\[.2pc]
 1&& \quad \quad 0&&0 \end{array}}\vline  \frac{\begin{array}{lcccr}
\sqrt{3}\quad \: \quad \quad&& 0 \quad \quad && 0\\ & & & &\\
\ & & & &\\ 0\quad \quad \quad \: && 0\quad \quad&&0 \\& & & &\\ & & & &\\0&& \quad
0\quad \quad \quad &&0 \end{array}}
{\begin{array}{lcccr}\: 0 \quad  &&0&& \quad   \sqrt{\frac{2}{\sqrt{3}}+1}\\
&&&& \\ \:  \sqrt{\frac{2}{\sqrt{3}}-1}  && \quad 0 &&0 \\ &&&& \\
0 & & \quad \sqrt{\frac{2}{\sqrt{3}}} && 0 \end{array}}\!\!\end{bmatrix}\!.
\label{eq7}
\end{equation}
It can be seen that $N$ is normal and plainly $\|B\| = \sqrt{3}$ while $\|C\|=1.$ When
$n=4$, it is impossible to find such a matrix, and that is our first theorem.

The following elementary lemma (which can be verified by induction on the integer $k$)
is used repeatedly in our proof.

\begin{lemm}
{\it Let $V$ be an $n$-dimensional vector space and let $V_1, \ldots, V_k$ be
subspaces of $V$ the sum of whose dimensions is larger than $(k-1)n;$ i.e.$,$
\begin{equation*}
\sum_{j=1}^k \dim \,  V_j \,  > \,  (k-1) n.
\end{equation*}
Then the intersection of these $k$ subspaces is nonzero.}
\end{lemm}

\begin{theor}[\!]
There exists a normal matrix $N$ of the form $\eqref{eq1}$ with
\begin{equation}
\|B\| = \sqrt{n} \ \|C\| \label{eq8}
\end{equation}
if and only if $n \leq 3.$
\end{theor}

\begin{proof}
Note first that if equalities \eqref{eq2} and \eqref{eq8} hold simultaneously, then
rank $B$ must be one and $C$ must be unitary. So, after applying a unitary similarity
by $\left[\begin{smallmatrix} C &O\\[.2pc] O &I \end{smallmatrix}\right],$ we may assume
that
\begin{equation}
N = \begin{bmatrix} A &B\\[.2pc]
I &D  \end{bmatrix}.\label{eq9}
\end{equation}
The normality condition $N^{*} N = NN^{*}$ leads to two equations
\begin{align}
A - D &= A^{*} B - BD^{*}, \label{eq10}\\[.3pc]
2I &= AA^{*} - A^{*}A + BB^{*} + B^{*}B+D^{*}D-DD^{*}.\label{eq11}
\end{align}
Since $B$ is of rank one,

where $\dim \,X$ stands for the dimension of a space $X.$ So, if $n \ge 3,$ then the
dimensions of $\ker \, B$ and $\ker \,B^{*}$ add up to more than $n.$ Hence their
intersection is nonzero, and we may choose a unit vector $x$ in this intersection. For
this vector, we obtain from \eqref{eq10}
\begin{equation}
(A-D) x =  - BD^{*} x, \label{eq12}
\end{equation}
and
\begin{equation}
(A-D)^{*} x = B^{*} Ax. \label{eq13}
\end{equation}
Equation \eqref{eq11} leads to the condition
\begin{equation}
2 = \|A^{*} x \|^2 - \|Ax \|^2 + \|Dx \|^2 - \|D^{*} x \|^2.\label{eq14}
\end{equation}

The rest of the proof shows that if $n > 3,$ then we can choose a vector $x \in (\ker
\, B) \cap (\ker \, B^{*})$ for which these conditions cannot be satisfied.

The two matrices $BD^{*}$ and $B^{*} A$ have rank at most 1, so their kernels have
dimension at least $n-1.$ Hence
\begin{equation}
\dim (\ker \, B) + \dim (\ker \, B^{*}) + \dim (\ker \, BD^{*}) + \dim (\ker \, B^{*}
A) \ge 4n - 4. \label{eq15}
\end{equation}
This is larger than $3n$ whenever $n > 4.$ So, in this case the four kernel spaces
involved in \eqref{eq15} have a nonzero intersection. Let $x$ be a unit vector in this
intersection. Then from \eqref{eq12} and \eqref{eq13} we find that
\begin{equation*}
(A-D)x = 0 \quad \mbox{and} \quad (A-D)^{*} x = 0.
\end{equation*}
Hence, $\|Ax\| = \|Dx \|$ and $\|A^{*}x \| = \|D^{*}x \|.$ This contradicts the
condition \eqref{eq14}.

Now consider the case $n=4.$ The spaces $\ker \, B$ and $\ker \, B^{*}$ have dimension
3 each, while the space $\ker \, B (A+D)^{*}$ has dimension at least 3. The three
dimensions add up to more than 8. Hence, we can find a unit vector $x$ in the
intersection of these three spaces. For this vector we have
\begin{align}
\|A^{*} x \|^2 - \|D^{*} x \|^2 &= \Re\, \langle (A+D)^{*} x, (A-D)^{*} x\rangle
\nonumber\\[.3pc]
&= \Re \, \langle (A+D)^{*}x, B^{*} Ax \rangle \nonumber\\[.3pc]
&= \Re \, \langle B (A+D)^{*} x, A x \rangle \nonumber\\[.3pc]
&= 0. \label{eq16}
\end{align}
Here the second equality is a consequence of \eqref{eq13}, and at the last step we
have used the fact that $B (A+D)^{*} x = 0.$

Using \eqref{eq12} instead of \eqref{eq13} we get
\begin{align}
\|Dx\|^2 - \|Ax\|^2 &=\Re \, \langle (A+D)x, (D-A) x \rangle \nonumber\\[.3pc]
&=\Re \, \langle (A+D) x, \textit{BD}^{*} x \rangle \nonumber\\[.3pc]
&=\Re \, \langle B^{*} (A+D)x, D^{*} x \rangle. \label{eq17}
\end{align}
Since $B$ is a matrix with rank equal to 1 and norm equal to 2, we
have $B^{*} BB^{*} = 4 B^{*}.$ (Use the polar decomposition
$B=\textit{UP}.$ In some orthonormal basis $P$ is diagonal with
only one nonzero entry 2 on the diagonal. So $B^{*} BB^{*} = P^3
U^{*} = 4 \textit{PU}^{*} = 4B^{*}.$) Hence we have
\begin{align*}
4 B^{*} A x &=B^{*} B B^{*} Ax\\[.3pc]
 &= B^{*} B (A-D)^{*}x \quad \mbox{(using}\, \eqref{eq13})\\[.3pc]
 &=B^{*} B (A+D)^{*} x - 2 B^{*} B D^{*} x\\[.3pc]
 &=-2 B^{*} B D^{*} x\\[.3pc]
 &=2 B^{*} (A-D)x \quad \mbox{(using}\, \eqref{eq12})\\[.3pc]
 &=4 B^{*} A x - 2 B^{*} (A+D)x.
\end{align*}
This shows that $B^{*} (A+D)x = 0,$ and we get from \eqref{eq17}
\begin{equation}
\|Dx\|^2 - \|Ax \|^2 = 0. \label{eq18}
\end{equation}
Clearly the relations \eqref{eq14}, \eqref{eq16} and \eqref{eq18} cannot be
simultaneously true.

We have shown that when $n \ge 4,$ there cannot exist a $2 n \times 2n$ normal matrix
of the form \eqref{eq9} in which $B$ is an $n \times n$ matrix of rank one. This
proves the theorem. \hfill $\qed$
\end{proof}

Our discussion leads to some natural questions.

\setcounter{theor}{0}
\begin{prob}
{\rm For $n \geq 4$, evaluate the quantity
\begin{equation*}
\alpha_n = \sup \left \{\|B\| / \|C\|\!\!: \exists A, D \
\mbox{for which}\ \left [\begin{array}{lr} A&B\\C&D
\end{array}\right ]\ \mbox{is normal} \right \}.
\end{equation*}
We have seen $\alpha_n < \sqrt{n}$ for $n \geq 4.$ It would be of interest to know
whether ${\alpha_n}$ is a bounded sequence.}
\end{prob}

\begin{prob}
{\rm What matrix pairs $B,C$ can be the off-diagonal entries of a normal matrix $N$ as
in \eqref{eq1}? In other words, when does $\left [\begin{smallmatrix} ? &B\\[.2pc] C &?
\end{smallmatrix}\right]$ have a normal completion?}
\end{prob}

\setcounter{theor}{0}
\begin{exam}
{\rm Consider the $2 \times 2$ matrices
\begin{equation*}
B = \begin{bmatrix} 1 &\varepsilon\\[.2pc] 0  &0 \end{bmatrix}, \quad C = \begin{bmatrix} 1 &0\\[.2pc]
0 & \varepsilon \end{bmatrix}.
\end{equation*}
Then, $\|B \|_2 = \| C\|_2.$ However, there do not exist any $2 \times 2$ matrices $A$
and $D$ for which $\left[ \begin{smallmatrix} A &B \\[.2pc] C& D \end{smallmatrix} \right]$
is normal. We leave the verification of this statement to the
reader. \ Thus the equality \eqref{eq2} is only a necessary
condition for normality of the matrix \eqref{eq1}.

We consider some special cases of the question raised in Problem~2. We assume either
$B=C,$ or $B=C^{*}.$

For every $B,$ the matrix $\left[ \begin{smallmatrix} ? &B\\[.2pc] B &? \end{smallmatrix}
\right]$ has a normal completion, and this completion may be chosen to be of the
special type $\left[ \begin{smallmatrix} A &B\\[.2pc] B &A \end{smallmatrix} \right].$
Indeed, if $U$ is the unitary matrix $U = \frac{1}{\sqrt{2}} \left[\begin{smallmatrix}
I &I\\[.2pc]
-I &I
\end{smallmatrix} \right],$ then
\begin{equation*}
U \begin{bmatrix}
A &B\\[.2pc]
B & A \end{bmatrix} U^{*} = \begin{bmatrix}
A+B &0\\[.2pc]
0 &A-B\end{bmatrix}.
\end{equation*}
So $\left[\begin{smallmatrix}
A &B\\[.2pc]
B &A \end{smallmatrix}\right]$ is normal if and only if $\left[ \begin{smallmatrix}
A+B &0\\[.2pc]
0 &A-B\end{smallmatrix}\right]$ is normal, and this is the case if and only if $A+B$
and $A-B$ both are normal. The most obvious choice of $A$ that assures this is
$A=B^{*}.$ Thus
\begin{equation}
\widetilde{B} = \begin{bmatrix}
B^{*} &B\\[.2pc]
B &B^{*}\end{bmatrix}\label{eq19}
\end{equation}
is a normal completion of $\left[ \begin{smallmatrix} ? &B\\[.2pc]
B &? \end{smallmatrix}\right].$ We have the norm inequality
\begin{equation}
\|B\| \leq \| \widetilde{B} \| \leq 2 \|B\|.\label{eq20}
\end{equation}
When $B = \left[ \begin{smallmatrix}
0 &1\\[.2pc]
0 &0\end{smallmatrix} \right]$ we have $\| \widetilde{B} \| = \|B\|.$ On the other
hand, if $B$ is any hermitian matrix, then $\|\widetilde{B} \| = 2 \|B\|.$ In this
case, and
more generally when $B$  is normal, $\left[ \begin{smallmatrix} 0 &B\\[.2pc]
B &0
\end{smallmatrix} \right]$ is normal and has norm equal to $\|B\|.$ This raises the
question of finding completions of $\left[
\begin{smallmatrix}
? &B\\[.2pc]
B & ? \end{smallmatrix}\right]$ that are `optimal' in various senses.}
\end{exam}

\setcounter{theor}{2}
\begin{prob}
{\rm Given a matrix $B$ find a matrix $A$ such that
\begin{equation*}
N = \begin{bmatrix} A &B \\[.2pc] B &A \end{bmatrix}
\end{equation*}
is normal and has the least possible norm. This is equivalent to
asking for a matrix $A$ such that $A+B$ and $A-B$ are normal and
the quantity $\max (\|A + B \|, \|A-B\|)$ is minimised. It might
be difficult to find {\it all} solutions to this problem. The
following considerations lead to {\it one} solution.

We assume that $B$ is a contraction, i.e. $\| B \| \leq 1$ and ask for an $A$ so that
$\left [ \begin{smallmatrix} A &B \\[.2pc] B &A \end{smallmatrix} \right]$ is unitary. This
is a unitary completion of the matrix $\left[ \begin{smallmatrix} ? &B\\[.2pc] B &?
\end{smallmatrix} \right].$ Let $B = \textit{US}\,\!V$ be the singular value decomposition of $B.$ Then
\begin{equation*}
\begin{bmatrix} U^{*} &0 \\[.2pc] 0 & U^{*} \end{bmatrix}  \begin{bmatrix} A &B \\[.2pc]
B &A \end{bmatrix} \begin{bmatrix} V^{*} &0 \\[.2pc] 0 &V^{*} \end{bmatrix}
= \begin{bmatrix} U^{*}AV^{*} &S\\[.2pc] S &U^{*}AV^{*} \end{bmatrix}.
\end{equation*}
So, the problem reduces to finding an $A^{\prime}$ such that
$\left[ \begin{smallmatrix} A^{\prime} & S \\[.2pc] S & A^{\prime} \end{smallmatrix} \right]$
is unitary. A familiar idea from the theory of unitary dilations
(p.~232 of \cite{2}) suggests the choice $A^{\prime} = i (I -
S^2)^{1/2}.$

This tells us how to find for any matrix $B$ one of the least-norm normal completions
of  $\left[ \begin{smallmatrix} ? &B\\[.2pc] B &? \end{smallmatrix} \right].$ Assume
$\|B\| = 1$ and find a unitary completion as proposed above.

Next we consider the case $B=C^{*},$ and ask for matrices $A$ and $D$ such that
\begin{equation}
N =  \begin{bmatrix} A &B\\[.2pc] B^{*} &D \end{bmatrix}\label{eq21}
\end{equation}
is normal. A calculation shows that the matrices $A$ and $D$ must be normal and
satisfy the equation
\begin{equation}
(A - A^{*})B = B (D - D^{*}).\label{eq22}
\end{equation}
Let $A = H_1 + i K_1$ and  $D = H_2 + i K_2$ be the Cartesian
decompositions of $A$ and $D.$ Here  $(H_1, K_1)$ and $(H_2, K_2)$
are two pairs of commuting hermitian matrices.
Equation~\eqref{eq22} is equivalent to $K_1 B = BK_2.$ This shows
that
\begin{equation*}
B^{*} B K_2 = B^{*} K_1 B = (K_1 B)^{*} B =  (BK_2)^{*}B = K_2
B^{*}B.
\end{equation*}
So $K_2$ commutes with $B^{*} B,$ and hence with the factor $P$ in
the polar decomposition $B=\textit{UP}.$

Thus the general solution to \eqref{eq22} is obtained as follows: Choose $K_0$ and
$K_2,$ both hermitian, satisfying the conditions
\begin{equation*}
K_0  P = \textit{PK}_0, \quad K_2 P = \textit{PK}_2, \quad (K_0 -
K_2) P = 0.
\end{equation*}
Let $K_1 = \textit{UK}_0 U^{*}.$ This condition ensures
\begin{equation*}
K_1 B = \textit{UK}_0U^{*}\!B = \textit{UK}_0P = \textit{UK}_2 P =
\textit{UPK}_2 = \textit{BK}_2.
\end{equation*}
Choose hermitian matrices $H_1$ and $H_2$ that commute with $K_1$ and $K_2,$
respectively. Let $A = H_1 + iK_1$ and $D = H_2 + iK_2.$ This leads to $N$ in
\eqref{eq21} being normal.

As before, we also consider the special case $\|B\| \leq 1$ and
ask for $A$ and $D$ such that the matrix \eqref{eq21} is unitary.
This can be solved as follows: Let $B = \textit{UP}$ be any polar
decomposition. Choose hermitian matrices $K_0$ and $K_2$ that
commute with $P$ and satisfy the inequalities
\begin{equation*}
K_0^2 \leq I - P^2, \quad K_2^2 \leq I - P^2.
\end{equation*}
Then choose hermitian matrices $H_0$ and $H_2$ that commute with $K_0$ and $K_2,$
respectively, and satisfy the conditions
\begin{equation*}
H_0^2 + K_0^2 = H_2^2 + K_2^2 = I - P^2.
\end{equation*}
Let $A = U (H_0 + i K_0) U^{*}$ and $D = H_2 + iK_2.$ Then the matrix \eqref{eq21} is
unitary.

Example~1 shows that the equality $\|B\|_2 = \|C\|_2$ is not a sufficient condition
for the existence of a normal completion of $\left[ \begin{smallmatrix} ? &B \\[.2pc] C &?
\end{smallmatrix} \right].$

Our next proposition shows that equality between all unitarily invariant norms is a
sufficient condition.}
\end{prob}

\begin{prop}$\left.\right.$\vspace{.5pc}

\noindent {\it Let $B, C$ be $n \times n$ matrices with
$|||B|||=|||C|||$ \ for every
unitarily invariant norm. Then the matrix $\left [\begin{smallmatrix} ? &B \\[.2pc]
C &? \end{smallmatrix}\right]$ has a completion that is a scalar multiple of a unitary
matrix.}
\end{prop}

\begin{proof}
If $|||B|||=|||C|||$ for every unitarily invariant norm, then $s_j (B) = s_j(C)$ for
all $j=1,2,\ldots, n.$ Hence, there exist unitary matrices $U_1, U_2, V_1, V_2$ such
that $B=U_1 SU_2,$ and $C=V_1 SV_2.$ Divide $B$ and $C$ by $||S||,$ and thus assume
$||S||=1.$ Then $I-S^2$ is positive, and has a positive square root. It is easy to see
that the matrix
\begin{equation*}
\begin{bmatrix} (I-S^2)^{\frac{1}{2}} &S \\[.2pc] S &-(I-S^2)^{\frac{1}{2}} \end{bmatrix}
\end{equation*}
is unitary. Multiply this matrix on the left by the unitary matrix $U_1 \oplus V_1,$ and on the right by the unitary matrix $V_2 \oplus U_2.$ This gives a unitary matrix whose off-diagonal blocks are $B$ and $C.$ \hfill{$\qed$}
\end{proof}

While the condition in the Proposition is not necessary, it is sensitive to small
perturbations. The matrices $B$ and $C$ in Example 1 satisfy the conditions $\|B\|_2 =
\|C\|_2,$ $|||B|||= |||C||| + O(\varepsilon),$ but for $\varepsilon \neq 0,$ there is
no possible normal completion of $\left[ \begin{smallmatrix} ? &B\\[.2pc] C &?
\end{smallmatrix} \right].$

\section*{Acknowledgement}

The second author thanks the Indian Statistical Institute and NSERC of Canada for
supporting a visit to New Delhi during which this work was initiated.

\end{document}